\documentclass{article}

\begin{document}

\title{Criticality in unbounded-types branching processes}

\author{G.T. Tetzlaff}

\date{}

\maketitle

Departamento de Computaci\'{o}n, FCEyN, Universidad de Buenos Aires

Ciudad Universitaria, 1428 Buenos Aires, Argentina

\begin{abstract}
Conditions for almost sure extinction are studied in discrete time branching processes with an infinite number of types. It is not assumed that the expected number of children is a bounded function of the parent's type. There might also be no integer $m$ such that there is a lower positive bound, uniform over the ancestor's type, for the probability that a population is extinct at the $m$-th generation. A weaker condition than the existence of such an $m$ is seen to lead to extinction almost surely if the sequence of expected generation sizes does not tend to infinity. Some criteria for a positive probability of nonextinction are given. Examples are provided by extending to our setting two applications, namely Leslie population dynamics and processes arising in continuum percolation in which the offsprings follow Poisson point distributions.
\end{abstract}



\section{Introduction}

We study the subject of giving conditions for the extinction probability to be
$1$ in discrete time branching processes with an infinite number of types.
Known results suppose that the expected number of children is a bounded
function of the type of their parent (Harris, 1963; Mode, 1971). Here we do
not prevent types from being arbitrarily large. This means that there might be
no uniform bound over the types for the expected number of children, nor a
number $m$\ of generations such that there is a uniform positive lower bound
for the probability that an individual will have no descent at the $m$-th
generation. A weaker condition than the latter is introduced in Section 2 and
it is shown to be sufficient for almost sure extinction if the sequence of
expected generation sizes does not tend to infinity, even if these
expectations fluctuate without an upper bound. Existence of a strictly
dominant eigenvalue of the integral operator associated to the offspring means
is not assumed in our propositions. Focusing on the expected number of
individuals in subsets of types that are in some sense larger than the type of
the initial individual, criteria for the extinction probability to lie below
$1$ are obtained in Section 3. Section 4 provides examples by extending to
unbounded types a couple of well known applications. The possibility of
constructing counterexamples shows that the sufficient condition for almost
sure extinction in Section 2 is not always necessary and leaves a criticality
problem open.

\section{Almost sure extinction}

Let $Z_{0},Z_{1},...$ be a discrete-time branching process with types in a
space $X$, countable or not. In the uncountable case, $X$ is an Euclidean
space and $Z_{0},Z_{1},...$ a sequence of point distributions that follows the
branching process definition as given in Harris (1963), Chapter 3. The number
of individuals with type in a measurable set $A\subset X$ at the $n$-th
generation will be denoted by $Z_{n}(A)$. The symbols $P_{x}()$, $E_{x}$ and
$q_{x}$ will stand for probabilities, expectations and the extinction
probability if the process is initiated with a single individual of type $x\in
X$.

\bigskip

\textbf{Condition 1} \ For any number of children $k>0$, there exist a number
of generations $m(k)$ and a real number $q(k)>0$ such that%
\[
\inf_{x\in X}P_{x}(Z_{m(k)}(X)=0|1\leq Z_{1}(X)\leq k)\geq q(k).
\]

\bigskip

\textbf{Proposition 2} \ If Condition 1 holds and $\lim\inf_{n}E_{x}%
Z_{n}(X)<\infty$, then $q_{x}=1$.

\bigskip

\textbf{Proof}. By Fatou's Lemma, $\lim\inf_{n}E_{x}Z_{n}(X)<\infty$ implies
\\$E_{x}\lim\inf_{n}Z_{n}(X)<\infty$, i.e. $\sum_{k\geq0}P_{x}(\lim\inf_{n}%
Z_{n}(X)>k)<\infty$. Thus given $\epsilon>0$, there is a positive integer
$k_{\epsilon}$ such that $P_{x}(\lim\inf_{n}Z_{n}(X)>k_{\epsilon})\leq
\epsilon$.

So $P_{x}(Z_{n}(X)>k_{\epsilon}$ for all but finite $n$'s$)\leq\epsilon$ and
therefore
\[
P_{x}(Z_{n}(X)\leq k_{\epsilon},i.o.)>1-\epsilon.
\]
We shall see that for any positive integer $k$, $P_{x}(1\leq Z_{n}(X)\leq
k,i.o.)=0$. This will complete the proof because then for all $\epsilon$,
$P_{x}(Z_{n}(X)=0,i.o.)>1-\epsilon$ and the extinction probability must equal
$1$. Let
\[
N_{0}(k)=0
\]
and for $i=1,2,...$,
\[
N_{i}(k)=\min\{n>N_{i-1}(k)-1+m(k):1\leq Z_{n}(X)\leq k\}
\]
with $\min\emptyset=\infty$. Now,
\[
P_{x}(1\leq Z_{n}(X)\leq k,i.o.)=P_{x}(\cap_{i}\{N_{i}(k)<\infty
\})=\lim_{i\rightarrow\infty}P_{x}(N_{i}(k)<\infty),
\]
so it is enough to prove that $P_{x}(N_{i}(k)<\infty)\leq(1-q(k)^{k})^{i}$.

For $i=0$ we have
\[
P_{x}(N_{0}(k)<\infty)=1=(1-q(k)^{k})^{0}%
\]
and for $i>0$%
\[
P_{x}(N_{i+1}(k)<\infty)\leq P_{x}(N_{i}(k)<\infty,Z_{N_{i}(k)-1+m(k)}(X)>0)
\]%
\[
=\sum_{n=1}^{\infty}P_{x}(Z_{n-1+m(k)}(X)>0|N_{i}(k)=n)P_{x}(N_{i}(k)=n).
\]
But given $N_{i}(k)=n$, there are at most $k$ individuals at generation $n$.
These must be the children of at most $k$ individuals out of generation $n-1$,
which have at most $k$ children each. Hence by Condition 1, $1-q(k)^{k}$ is an
upper bound for $P_{x}(Z_{n-1+m(k)}(X)>0|N_{i}(k)=n)$ yielding
\[
P_{x}(N_{i+1}(k)<\infty)\leq(1-q(k)^{k})\sum_{n=1}^{\infty}P_{x}(N_{i}(k)=n)
\]%
\[
=(1-q(k)^{k})P_{x}(N_{i}(k)<\infty).
\]
The desired result follows by induction.

\section{Positive probability of nonextinction}

Let $Z_{0},Z_{1},...$ be a discrete-time branching process with types in a
countable or uncountable set $X$ as in the preceding section. Given a type $x$
and a positive integer $n$, by the $n$-th generation descendants of $x$ we
will mean the individuals in the $n$-th generation of a population that begins
with a single individual of type $x$. Let us consider the following order
relation in the set of types. Given two types $x$ and $y$ we shall say that
$x$ is larger than $y$ if for every measurable subset $A$ of types and any $n$
$(n=1,2,...)$, the number of $n$-th generation descendants with types in $A$
is stochastically larger for $x$ than for $y$.

\bigskip

\textbf{Proposition 3} \ Let $y$ be a type in $X$. Suppose that there exist a
measurable set $A$ whose elements are larger than $y$ and a positive integer
$i$ such that $E_{y}Z_{i}(A)>1$. Then for any $x$ larger than $y$ we have
$q_{x}<1$.

\bigskip

\textbf{Proof.} Let $S_{0},S_{1},...$ be a single-type branching process
defined by
\[
S_{0}=1
\]
and for $k=0,1,...$,
\[
P(S_{1}=k)=P_{y}(Z_{i}(A)=k).
\]
As $E_{y}Z_{i}(A)>1$, this process is supercritical. Now consider $Z_{0}%
,Z_{i},Z_{2i},...$, the multitype process at those generations that are
multiples of $i$, beginning with a single individual whose type is $x$. Since
$x$ and all elements of $A$ are larger than $y$, $Z_{ni}(A)$ is stochastically
larger than $S_{n}$ for any $n$. As $Z_{ni}(X)\geq Z_{ni}(A)$, we have
$q_{x}<1$.

\bigskip

\textbf{Corollary 4} \ If for a type $x$ there exists $i$ such that
$E_{x}Z_{i}(\{x\})>1$ then $q_{x}<1$.

\bigskip

\textbf{Corollary 5} \ Suppose that there exist $y$, $i$ and $A$ as in the
proposition. If there is a type $x$, an integer $j$ and a set $B\subset X$
whose elements are larger than $y$, such that $E_{x}Z_{j}(B)>0$, then
$q_{x}<1$.

\bigskip

\textbf{Corollary 6} \ If there is a smallest type $s$ in the sense that every
type $x$ is larger than $s$, and $E_{s}Z_{i}(X)>1$ for some $i$, then
$q_{x}<1$ for every $x$.

\section{Examples}

Before going into the examples it should be remarked that in all the processes
we are dealing with, the expectations $E_{x}Z_{n}(A)$ verify%
\[
E_{x}Z_{n+1}(A)=\int_{X}E_{y}Z_{n}(A)dE_{x}Z_{1}(y),
\]
for $n=1,2,...$, as proved in Harris (1963), Chapter 3, and in Mode (1971),
Chapter 6. The set $A$ is any measurable subset of $X$ and $dE_{x}Z_{1}(y)$
means that the integration is performed with respect to the measure defined by
$E_{x}Z_{1}(S)$ for any measurable subset $S$ of $X$. Under discrete
interpretation, the right hand side of the equality becomes the result of
rising the first generation mean matrix to the $(n+1)$-th power and taking the
sum of the elements in the row of $x$ that lie in the columns that correspond
to types in $A$. In the continuous case, suppose that $k(x,y)$ is a nucleus
such that for any measurable $A\subset X$ we have $E_{x}Z_{1}(A)=\int
_{A}k(x,y)d\mu(y)$ for some measure $\mu$ on $X$. Then the formula yields the
iteration of the integral operator $K$ associated to $k$ and $\mu$, allowing
us to write $E_{x}Z_{n+1}(A)=\int_{A}K^{n}(k(x,y))d\mu(y)$.

\bigskip

In all our examples we shall both have%
\[
\sup_{x\in X}E_{x}Z_{1}(X)=\infty
\]
as%
\[
\inf_{x\in X}P_{x}(Z_{m}(X)=0)=0
\]
for any $m$. That the expectations are not bounded will be evident. In order
to prove the second property, the following proposition will be useful.

\bigskip

\textbf{Proposition 7} \ Suppose that for any positive integer $h$,
\\$\inf_{x\in X}P_{x}(Z_{1}(X)\leq h)=0$ and that there exists $\alpha>0$ such
that for every $x\in X$, $P_{x}(Z_{1}(X)>0)\geq\alpha$. Then $\inf_{x\in
X}P_{x}(Z_{m}(X)=0)=0$ for any number of generations $m$.

\bigskip

\textbf{Proof.} Let $\{x_{h}:h=1,2,...\}\subset X$ be such that $P_{x_{h}%
}(Z_{1}(X)\leq h)\rightarrow0$ when $h\rightarrow\infty$. For any $m$,%
\[
P_{x_{h}}(Z_{m}(X)=0)
\]
\[
\leq P_{x_{h}}(Z_{1}(X)\leq h)+P_{x_{h}}(Z_{m}(X)=0|Z_{1}(X)>h)P_{x_{h}}%
(Z_{1}(X)>h)
\]%
\[
\leq P_{x_{h}}(Z_{1}(X)\leq h)+(1-\alpha^{m-1})^{h}.
\]
So $\lim\inf_{h}P_{x_{h}}(Z_{m}(X)=0)\leq0$ and the proposition is proved.

\subsection{Leslie dynamics}

We shall consider the infinite types analog of the finite-types branching
process studied by Pollard (1966), which is the stochastic version of the
deterministic population dynamics described by Leslie (1945). The set of types
will be here $\{x_{0},x_{1},...\}$. If we adopt the interpretation of types as
age classes, then the individuals of type $x_{0}$ are the youngest ones. An
$x_{0}$ may bear individuals of type $x_{0}$ or may produce an individual of
type $x_{1}$ by surviving enough time to enter in the next age class. In
general $x_{i}$ $(i=0,1,...)$ may only produce $x_{0}$'s and eventually an
$x_{i+1}$, so that the infinite mean matrix has only two positive entries at
each row, $E_{x_{i}}Z_{1}(\{x_{0}\})$ and $E_{x_{i}}Z_{1}(\{x_{i+1}\})$. For
simplicity, both the individuals $x_{0}$ as the $x_{i+1}$ arising from an
$x_{i}$ will be called its children.

\bigskip

Our examples will have first generation means $E_{x_{i}}Z_{1}(\{x_{0}%
\})=p^{1-i}$ and $E_{x_{i}}Z_{1}(\{x_{i+1}\})=p^{2}$ for $p$ in $(0,1)$, so
the infinite mean matrix is%
\[
\left(
\begin{array}
[c]{cccccc}%
p & p^{2} & 0 & 0 & 0 & ...\\
1 & 0 & p^{2} & 0 & 0 & \\
p^{-1} & 0 & 0 & p^{2} & 0 & \\
p^{-2} & 0 & 0 & 0 & p^{2} & \\
... &  &  &  &  & ...
\end{array}
\right)  .
\]
By what we noticed at the beginning of this section, the first row of the
$n$-th power of this matrix consists of the $n$-th generation means $E_{x_{0}%
}Z_{n}(\{x_{0}\})$, $E_{x_{0}}Z_{n}(\{x_{1}\})$,.... It can be seen by
induction over $n$ that these entries are
\[
E_{x_{0}}Z_{n}(\{x_{i}\})=\left\{
\begin{array}
[c]{ccc}%
(1/2)(2p)^{n}(p/2)^{i} & if & 0\leq i<n\\
p^{2n} & if & i=n\\
0 & if & i>n
\end{array}
\right.  .
\]
In particular, $E_{x_{0}}Z_{n}(\{x_{0}\})=(1/2)(2p)^{n}$ so by Corollary 4,
$q_{x_{0}}<1$ if $p>1/2$. By a simple case of Corollary 5, $q_{x_{i}}<1$ for
any $i$.

\bigskip

In case $p\leq1/2$ the mean number of individuals of any type at the $n$-th
generation $(n=1,2,...)$ in a population that started with a single individual
of type $x_{0}$ is bounded by a constant. Indeed
\[
E_{x_{0}}Z_{n}(X)=\sum_{i=0}^{n}E_{x_{0}}Z_{n}(\{x_{i}\})=(1/2)(2p)^{n}%
\frac{1-(p/2)^{n}}{1-p/2}+p^{2n}<1.
\]
Thus $\lim\inf_{n}E_{x_{0}}Z_{n}(X)<\infty$ is satisfied. But\ if $p<1/2$ we
can bound the expected total number of descendants in all generations because
\[
\sum_{n=1}^{\infty}E_{x_{0}}Z_{n}(X)<\frac{1}{2-p}\sum_{n=1}^{\infty}%
(2p)^{n}+\sum_{n=1}^{\infty}p^{2n}<\infty.
\]
So $E_{x_{0}}\sum_{n=1}^{\infty}Z_{n}(X)<\infty$ and therefore $q_{x_{0}}=1$.

\bigskip

We shall now consider two different offspring distributions behind our means.
In both examples, the number of $x_{0}$-type children of an $x_{i}$ parent
$(i=0,1,...)$ will follow a Poisson distribution with mean $p^{1-i}$ and the
number of $x_{i+1}$-type children will be $1$ with probability $p^{2}$ and $0$
with probability $1-p^{2}$. The difference will lie in the joint distribution
for these marginals, as described below. It can already be seen that
$\inf_{x\in X}P_{x}(Z_{m}(X)=0)=0$ for any $m$ because having these marginal
distributions implies the hypotheses of Proposition 7.

\bigskip

The first choice of distributions for the children of the $x_{i}$'s
$(i=0,1,...)$ is one that satisfies Condition 1. We take for each $x_{i}$ the
joint distribution that makes the number of children of type $x_{0}$
independent of the number of $x_{i+1}$ children. This family satisfies
Condition 1 for $m(k)=2$ because%
\[
P_{x_{i}}(Z_{2}(X)=0|1\leq Z_{1}(X)\leq k)>
\]%
\[
>P_{x_{i}}(Z_{2}(X)=0|Z_{1}(\{x_{i+1}\})=0,1\leq Z_{1}(X)\leq k)\cdot
\]%
\[
\cdot P_{x_{i}}(Z_{1}(\{x_{i+1}\})=0|1\leq Z_{1}(X)\leq k)
\]
where%
\[
P_{x_{i}}(Z_{2}(X)=0|Z_{1}(\{x_{i+1}\})=0,1\leq Z_{1}(X)\leq k)\geq(P_{x_{0}%
}(Z_{1}(X)=0))^{k},
\]
while for any integer $h$ between $1$ and $k$,%
\[
P_{x_{i}}(Z_{1}(\{x_{i+1}\})=0|Z_{1}(X)=h)=
\]%
\[
=\frac{(1-p^{2})e^{-p^{1-i}}(p^{1-i})^{h}/h!}{(1-p^{2})e^{-p^{1-i}}%
(p^{1-i})^{h}/h!+p^{2}e^{-p^{1-i}}(p^{1-i})^{h-1}/(h-1)!}>\frac{1}{1+k}.
\]
Thus%
\[
\inf_{x_{i}\in X}P_{x_{i}}(Z_{2}(X)=0|1\leq Z_{1}(X)\leq k)>(P_{x_{0}}%
(Z_{1}(X)=0))^{k}\frac{1}{1+k}.
\]
Since it has been shown that $\lim\inf_{n}E_{x_{0}}Z_{n}(X)<\infty$ for
$p\leq1/2$, we have $q_{x_{0}}=1$ by Proposition 2.

\bigskip

The second family of distributions does not satisfy Condition 1 because there
is enough positive correlation between bearing few children and having one
very prolific among them. A way to achieve this is by requiring that at least
for $k=1$ and for infinite $i$'s,
\[
P_{x_{i}}(Z_{1}(\{x_{i+1}\})=1|1\leq Z_{1}(X)\leq k)=1.
\]
This is possible because due to the marginal distributions $P_{x_{i}}%
(Z_{1}(X)\leq1)\rightarrow0$ when $i\rightarrow\infty$. So for any $i$ greater
than some $\widetilde{i}$, $P_{x_{i}}(Z_{1}(X)\leq1)<p^{2}$. Since $p^{2}$ is
the probability of having one $x_{i+1}$, we can define for those $i$'s a joint
distribution that loads all of $P_{x_{i}}(Z_{1}(X)\leq1)$ only on the
probability of having no $x_{0}$ and one $x_{i+1}$. So for $i>\widetilde{i}$,
$P_{x_{i}}(Z_{1}(\{x_{i+1}\})=1|Z_{1}(X)=1)=1$, as we required. Now we are
satisfying an analogous hypothesis as that of Proposition 7, namely that for
any integer $h$, $\inf_{x_{i}\in X}P_{x_{i}}(Z_{2}(X)\leq h|Z_{1}(X)=1)=0$. By
this and the fact that there exists $\alpha>0$ such that for every $x_{i}$,
$P_{x_{i}}(Z_{1}(X)>0)>\alpha$, an analogous proof as that of Proposition 7
gives%
\[
\inf_{x_{i}\in X}P_{x_{i}}(Z_{m}(X)=0|Z_{1}(X)=1)=0
\]
for any number of generations $m$. This means that Condition 1 is not
satisfied because it fails for $k=1$. Since it was already proved that
$q_{x_{0}}=1$ for $p<1/2$, a consequence of this example is that the
hypotheses of Proposition 2 are not necessary for the extinction probability
to be $1$. For $p=1/2$, the criticality problem remains here unsolved.

\bigskip

Finally it can be said about the case $p\leq1/2$ for $i=1,2,...$, that
$q_{x_{i}}=1$ if $q_{x_{0}}=1$ by types communication. In fact, $q_{x_{i}}<1$
implies $q_{x_{i-1}}<1$.

\subsection{A non recurrent variant}

The following modification of the unbounded-types Leslie scheme provides an
example of the application of Proposition 2 in a non recurrent case. Let
$X=\{x_{0},x_{1},...\}$, $p\in(0,1)$ and the infinite mean matrix be%
\[
\left(
\begin{array}
[c]{cccccc}%
p & 0 & 0 & 0 & 0 & ...\\
1 & 0 & p & 0 & 0 & \\
p^{-1} & 0 & 0 & p & 0 & \\
p^{-2} & 0 & 0 & 0 & p & \\
... &  &  &  &  & ...
\end{array}
\right)  ,
\]
i.e. the only positive expectations are $E_{x_{i}}Z_{1}(\{x_{0}\})=p^{1-i}$
for $i=0,1,...$ and $E_{x_{i}}Z_{1}(\{x_{i+1}\})=p$ for $i=1,2,...$. In a
population that begins with a single individual of type $x_{1}$, only types
$x_{0}$ and $x_{n+1}$ can be present at the $n$-th generation and it can
easily be seen that%
\[
E_{x_{1}}Z_{n}(\{x_{0}\})=\sum_{i=0}^{n-1}p^{i}%
\]
and%
\[
E_{x_{1}}Z_{n}(\{x_{n+1}\})=p^{n}.
\]
Hence%
\[
E_{x_{1}}Z_{n}(X)=\sum_{i=0}^{n}p^{i}<\frac{1}{1-p}.
\]

In the same way as in the preceding examples, we can find a family of
distributions with the given means that make $q_{x_{1}}=1$ by Proposition 2,
and another family of distributions can be defined, such that the criticality
problem is not solved by Proposition 2. Anyway we notice that although the
line of descent $x_{i},x_{i+1},...$ $(i\geq1)$ is a source of growing
quantities of individuals of type $x_{0}$, it is clearly once broken with
probability $1$. So as $q_{x_{0}}=1$ we have $q_{x_{i}}=1$ for any $i\geq1$.

\subsection{Poisson process descent}

Continuum percolation clusters in ${R}^{d}$ can be dominated by
coupling with branching processes so that extinction in the branching process
implies that the percolation cluster is finite. This is done by defining a
parent-child relation between objects that intersect each other in the
percolation setting. Extra individuals are eventually added to the families so
defined, so that the resulting families follow the independent offspring
distributions of a branching process. This method appears in Hall (1985) and
it is used and referenced in Meester and Roy (1996). It is a way to obtain an
upper bound for $\lambda$, the Poisson intensity parameter of the objects'
centers, such that below this bound the finiteness of the clusters is almost sure.

\bigskip

This application motivates us to analyze the criticality subject in branching
processes that would dominate continuum percolation clumps of objects whose
size might be arbitrarily large. We shall consider the case of spheres in
${R}^{d}$ ($d\geq1$) with not uniformly bounded random radii, which
will be the types in the branching process. Given a sphere in ${R}^{d}$,
its offspring is determined by spheres that have centers in a homogeneous
Poisson process in ${R}^{d}$ with intensity $\lambda$ and have iid
radii. The children are those spheres that intersect the given sphere. The
offspring of each of these children is generated by a new independent Poisson
process and independent radii an so on. It can be checked that we both have
$\sup_{x\in X}E_{x}Z_{1}(X)=\infty$ as $\inf_{x\in X}P_{x}(Z_{m}(X)=0)=0$ for
any $m$ by Proposition 7.

\bigskip

Let $F$ be the cumulative distribution function of the radii and assume that
it has an inverse $F^{-1}:[0,1)\rightarrow\lbrack0,\infty)$ when restricted to
$[0,\infty)$. A way to deduce the expression of a mean nucleus is to realize
children spheres as homogenous Poissonian points with intensity $\lambda$ in
${R}^{d}\times\lbrack0,1)$. The first $d$ dimensions give the center of
the child sphere and the last one is the image of its radius by $F$. Let
$u\in\lbrack0,1)$ be a value for the last component of a point in our
representation. We know that spheres with radius $F^{-1}(u)$ intersect a
parent sphere of radius $x$ if their centers lie inside the parent sphere or
outside of it but at an ${R}^{d}$ distance within $(x,x+F^{-1}(u)]$ to
the parent's center. So in our $(d+1)$-dimensional representation, the points
that correspond to $F^{-1}(u)$-type children of an $x$ lie in the
$d$-dimensional sphere that has the same center components as the center of
the parent sphere, radius $x+F^{-1}(u)$ and last component $u$. For any
$x_{1}\in(0,\infty)$, the union of such spheres over $u$ in the interval
$[F(0),F(x_{1})]$ is a solid with volume
\[
\int_{F(0)}^{F(x_{1})}\upsilon_{d}(x+F^{-1}(u))^{d}du=\int_{0}^{x_{1}}%
\upsilon_{d}(x+y)^{d}dF(y),
\]
where $\upsilon_{d}$ is the volume of the unit sphere in ${R}^{d}$.
Thus we obtained the volume of the subset of ${R}^{d+1}$ whose points
represent those children of an $x$ that have last component $u\in\lbrack
F(0),F(x_{1})]$ or equivalently, type in $[0,x_{1}]$. As our Poisson process
is homogeneous with parameter $\lambda$, the expected number of children of an
$x$ with types in $[0,x_{1}]$ is
\[
E_{x}Z_{1}([0,x_{1}])=\lambda\int_{0}^{x_{1}}\upsilon_{d}(x+y)^{d}dF(y).
\]

Let $k_{d}(x,y)=\lambda\upsilon_{d}(x+y)^{d}$. If we suppose that the
distribution of the radii has finite $2d$ moment (i.e. sphere content has
finite variance) then the integral operator $K_{d}:L^{2}([0,\infty
),dF)\rightarrow L^{2}([0,\infty),dF)$ defined by $k_{d}$ is compact because
$k_{d}$ is a Hilbert-Schmidt nucleus. Indeed, $\int_{0}^{\infty}\int
_{0}^{\infty}k_{d}^{2}(x,y)dF(x)dF(y)<\infty$. Since $K_{d}$ is also strictly
positive, it has a dominant eigenvalue $\rho>0$ and therefore by $E_{x}%
Z_{n}([0,\infty))=\int_{0}^{\infty}K_{d}^{n-1}(k_{d}(x,y))dF(y)$, the sequence
$\{E_{x}Z_{n}([0,\infty)):n=1,2,...\}$ is bounded if and only if $\rho\leq1$,
for any initial type $x$.

\bigskip

It will now be seen that for any dimension $d$ and every type $x$, we have
$\rho\leq1$ if and only if $q_{x}=1$. To complete the hypotheses of
Proposition 2 when $\rho\leq1$ it only remains to show that Condition 1 is
satisfied. Fixing $\widetilde{x}>0$, we shall first prove that there is a
positive bound for the probability of having all children with types below
$\widetilde{x}$, or more precisely that for any $k\geq1$,%
\[
\inf_{x\in\lbrack0,\infty)}P_{x}(Z_{1}([0,\infty))=Z_{1}([0,\widetilde
{x}])|1\leq Z_{1}([0,\infty))\leq k)>b(k)
\]
for some $b(k)>0$. By the representation of children as homogeneous Poissonian
points, we see that given a number of children, each choice of type for a
child is independent of the others and the probability that the choice falls
below $\widetilde{x}$ when the parent type is $x$, equals the volume ratio%

\[
\frac{\int_{0}^{\widetilde{x}}\upsilon_{d}(x+y)^{d}dF(y)}{\int_{0}^{\infty
}\upsilon_{d}(x+y)^{d}dF(y)}=\frac{x^{d}\int_{0}^{\widetilde{x}}%
dF(y)+\sum_{i=1}^{d}(_{i}^{d})x^{d-i}\int_{0}^{\widetilde{x}}y^{i}dF(y)}%
{x^{d}+\sum_{i=1}^{d}(_{i}^{d})x^{d-i}\int_{0}^{\infty}y^{i}dF(y)}.
\]
The integrals are finite by the assumption about the moments, so this
probability is a continuous and strictly positive function of $x$. Letting
$x\rightarrow\infty$ gives the strictly positive value $\int_{0}%
^{\widetilde{x}}dF(y)$. Hence the function has a positive lower bound
$\beta_{\widetilde{x}}$ and we can take $b(k)=\beta_{\widetilde{x}}^{k}$.
Finally, since an individual whose type is in $[0,\widetilde{x}]$ has no
children with a probability larger or equal than that for $\widetilde{x}$,
Condition 1 is satisfied by
\[
\inf_{x\in\lbrack0,\infty)}P_{x}(Z_{2}([0,\infty))=0|1\leq Z_{1}%
([0,\infty))\leq k)\geq\beta_{\widetilde{x}}^{k}(P_{\widetilde{x}}%
(Z_{1}(X)=0))^{k}.
\]

\bigskip

If $\rho>1$, taking $x=0$ as the initial type, $E_{0}Z_{i}([0,\infty))>1$ for
some generation $i$. By Corollary 6, we have $q_{x}<1$ for every $x\in
\lbrack0,\infty)$ if we prove that $0$ is a smallest type. In fact, $x\geq y$
implies that $x$ is larger than $y$ in the sense of Section 3. This assertion
can be justified by coupling, constructing a process that begins with an
individual of type $y$ jointly with one that begins with an $x$. We first
assign children to the $y$ by placing points with intensity $\lambda$ in the
solid of ${R}^{d+1}$ that corresponds to $y$'s first generation
descent, i.e. the union over $u\in\lbrack0,1)$ of the ${R}^{d}$ spheres
of radius $y+F^{-1}(u)$. This set is contained in the solid corresponding to
the first generation descent if the parent's type were $x$ because
$y+F^{-1}(u)\leq x+F^{-1}(u)$. So we can assign to $x$ the same children as
those of $y$ plus eventually some more in order to follow the children
distribution of an $x$ parent. In the subsequent generations, $x$ receives the
same descent as that of $y$ plus the descent of those first generation extra
children. This makes every $n$-th generation at any $A\subset X$, equal or
larger for $x$ than for $y$.

\bigskip

\textbf{Acknowledgement}

\bigskip

I am grateful to P.Ferrari for some useful advice.

\bigskip

\textbf{References}

\bigskip

Hall P., 1985. On continuum percolation. Ann. Probab. 13 4: 1250-1266.

Harris, T.E., 1963. The theory of branching processes. Spinger, Berlin.

Leslie, P.E., 1945. On the use of matrices in certain population mathematics.
Biometrika 33: 183-212.

Meester, R., Roy, R., 1996. Continuum percolation. Cambridge University Press.

Mode, C.J., 1971. Multitype branching processes: theory and applications.
American Elsevier Publishing Company, New York.

Pollard, J.H., 1966. On the use of the direct matrix product in analyzing
certain stochastic population models. Biometrika 53: 397-415.
\end{document}